\documentclass[12pt]{article}

\usepackage{amssymb}
\usepackage{amsmath}

\usepackage[latin1]{inputenc}
\usepackage[dvips]{epsfig}

\newtheorem{thm}{Theorem}[section]
\newtheorem{pro}[thm]{Proposition}
\newtheorem{lem}[thm]{Lemma}

\newtheorem{cor}[thm]{Corollary}
\newtheorem{rem}{Remark}

\newcommand{\SM}{{\mathcal S}}
\newcommand{\V}{{\mathcal V}}
\newcommand{\LE}{{\mathcal L}}
\newcommand{\la}{{\lambda}}

\newcommand{\Id}{\operatorname{Id}}
\newcommand{\Int}{\operatorname{Int}}
\newcommand{\Ad}{\operatorname{Ad}}
\newcommand{\For}{\operatorname{For}}
\newcommand{\NN}{\mathbb{N}}
\newcommand{\ZZ}{\mathbb{Z}}
\newcommand{\zero}{\hat{0}}
\newcommand{\one}{\hat{1}}
\newcommand{\rk}{\operatorname{rk}}
\newcommand{\join}{\vee}
\newcommand{\AAA}{\textsf{A}}
\newcommand{\BBB}{\textsf{B}}

\newenvironment{proof}{\begin{trivlist}\item{\bf{Proof.}}}
  {\hfill\rule{2mm}{2mm}\end{trivlist}}

\begin{document}
\title{Supersolvable LL-lattices of binary trees}
\author{Riccardo Biagioli and Fr\'ed\'eric Chapoton}
\date{\today}
\maketitle

\begin{abstract}
  Some posets of binary leaf-labeled trees are shown to be
  supersolvable lattices and explicit EL-labelings are given. Their
  characteristic polynomials are computed, recovering their known
  factorization in a different way.
\end{abstract}

\section{Introduction}

The aim of this article is to continue the study of some posets on
forests of binary leaf-labeled trees introduced by the second author
in \cite{bessel}. These posets have already been shown in \cite{fred}
to have nice properties. The main result there was the fact that the
characteristic polynomials of all intervals in these posets factorize
completely with positive integer roots. By a theorem of Stanley
\cite{staSS}, this property is true in general for the so-called
semimodular supersolvable lattices. Since these intervals are not
semimodular in general, one can not use this theorem to recover the
result of \cite{fred}. For a class of lattices, called LL-lattices,
containing the semimodular-supersolvable ones, a theorem due to Blass
and Sagan \cite{blassagan} generalizes Stanley's theorem.

The first main theorem of our article states that these intervals are
indeed lattices, which was not known before. The proof uses a new
description of the intervals using admissible partitions. Our second
main result is the fact that these lattices are supersolvable. We prove it
by giving explicit $S_n$ EL-labelings and using the recent criterion
of McNamara \cite{macnamara}. As third result, we show that these
intervals are LL-lattices and, using the theorem of Blass and Sagan mentioned above, we give a different proof of the factorization of
characteristic polynomials and the explicit description of roots which
were found in \cite{fred}.

\section{Notation, definitions and preliminaries}

In this section we give some definitions, notation and results that
will be used in the rest of this work. Let $\NN := \{1,2,3, \ldots \}$
and $\ZZ$ the set of integers. For every $n \in \NN$ let $[n]:=\{1,2,\ldots, n\}$. The
cardinality of a finite set $A$ is denoted by $|A|$.

\subsection{Posets}

We follow Chapter 3 of \cite{staEC1} for any undefined notation and
terminology concerning posets. Given a finite poset $(P,\leq)$ and
$x,y \in P$ with $x\leq y$ we let $[x,y]:=\{z \in P \; : \; x\leq z
\leq y \}$ and call this an {\em interval} of $P$. We denote by
$\Int(P)$ the set of all intervals of $P$. We say that $y$ {\em
  covers} $x$, denoted $x \lhd y$, if $|[x,y]|=2$. A poset is said to
be {\em bounded} if it has one minimal and one maximal element,
denoted by $\zero$ and $\hat{1}$ respectively. The {\em Möbius
  function} of $P$, $\mu: \Int(P) \rightarrow \ZZ$, is defined
recursively by

$$ \mu(x,y) := \left \{\begin{array}{ll}
1 & \mbox{if $x=y$}, \\
- \sum_{x \leq z < y} \mu(x,z) & \mbox{if $x \neq y$}.
\end{array} \right. $$
\smallskip

If $x,y \in P$ are such that $\{ z \in P : \; z \geq x, z \geq y \}$
has a minimum element then we call it the {\em join} of $x$ and $y$,
denoted by $x \join y$. Similarly, we define the {\em meet} of $x$ and
$y$ if $\{ z \in P: \; z \leq x , \; z \leq y\}$ has a maximum
element, denoted by $x \wedge y$. A {\em lattice} is a poset $L$ for
which every pair of elements has a meet and a join. A well-known
criterion is the following (see e.g. \cite[Proposition
3.3.1]{staEC1}).
\begin{pro}
  \label{lattice}
  If $P$ is a finite poset with $\hat{1}$ such that every pair of
  elements has a meet then $P$ is a lattice.
\end{pro}
A lattice $L$ that satisfies the following condition
\begin{equation}
  \label{semimodular}
  {\rm if} \;  x \; {\rm and} \; y \; {\rm both  \; cover} \; x \wedge y,\;{\rm then}  \; x \join y\;{\rm  covers \;both} \; x \; {\rm and} \; y,
\end{equation}
is said to be {\em semimodular}. The set of {\em atoms} of a finite lattice $L$, \textit{i.e.} the elements $a$ covering $\hat{0}$, is denoted by $\AAA(L)$. 

\subsection{Edge-labelings}

If $x,y \in P$, with $x \leq y$, a {\em chain} from $x$ to $y$ of {\em
  length} $k$ is a $(k+1)$-tuple $(x_0,x_1,\ldots,x_k)$ such that
$x=x_0<x_1<\ldots<x_k=y$. A chain $x_0\lhd x_1 \lhd \ldots \lhd x_k$
is said to be {\em saturated}. A poset $P$ with a $\zero$ is said to be
\textit{graded} if, for any $x \in P$, all saturated chains from
$\zero$ to $x$ have the same length, called the {\em rank} of $x$ and
denoted by $\rk(x)$. We denote by ${\mathcal M}(P)$ the set of all
maximal chains of $P$.

A function $\la:\{(x,y)\in P^2 : x \lhd y\} \rightarrow \NN$ is an
{\em edge-labeling} of $P$. For any saturated chain $m:x=x_0\lhd x_1
\lhd \ldots \lhd x_k=y$ of the interval $[x,y]$ we set
\[\la(m)=(\la(x_0,x_1),\la(x_1,x_2),\ldots,\la(x_{k-1},x_k)).\]
The chain $m$ is said to be {\em increasing} if $\la(x_0,x_1)\leq
\la(x_1,x_2)\leq \cdots \leq \la(x_{k-1},x_k)$. Let $\leq_L$ be the
lexicographic order on finite integer sequences, \textit{i.e.}
$(a_1,\ldots,a_k)<_L (b_1,\ldots, b_k)$ if and only if $a_i < b_i$
where $i=$ min$\{j \in [k]: a_j \neq b_j\}$.

An edge-labeling of $P$ is said to be an EL-$labeling$ if the
following two conditions are satisfied:
\begin{itemize}
\item[$i)$] Every interval $[x,y]$ has exactly one increasing saturated
chain $m$.
\item[$ii)$] Any other saturated chain $m^{\prime}$ from $x$ to $y$
satisfies $\la(m)<_L \la(m^{\prime})$.
\end{itemize}

A graded poset is said to be {\em edge-wise lexicographically
  shellable} or EL-$shellable$, if it has an EL-labeling.
EL-shellable posets were first introduced by Björner \cite{bjorner}.
Several connections with shellable, Cohen-Macaulay complexes and
Cohen-Macaulay posets can be found in the survey paper
\cite{bjgarsia}. In particular EL-shellable posets are Cohen-Macaulay
\cite{bjorner}.

\smallskip A particular class of EL-labelings has an interesting
property.

An EL-labeling $\la$ is said an $S_n$ EL-\textit{labeling} if, for any
maximal chain $m:\zero=x_0\lhd x_1 \lhd \cdots \lhd x_n=\hat{1}$ of
$P$, the label $\la(m)$ is a permutation of $[n]$. If a poset $P$ has
an $S_n$ EL-labeling, then it is said to be $S_n$
EL-\textit{shellable}.

Following \cite{staSS}, we introduce the following definition. A
finite lattice $L$ is said to be {\em supersolvable} if it contains a
maximal chain, called an {\em $M$-chain} of $L$, which together with any
other chain in $L$ generates a distributive sublattice. Examples of
supersolvable lattices include modular lattices, the partition lattice
$\Pi_n$ and the lattice of subgroups of a finite supersolvable group.

McNamara \cite[Theorem 1]{macnamara} has recently shown that
supersolvable lattices are completely characterized by $S_n$
EL-shellability.

\begin{thm}
  \label{MN}
  A finite graded lattice of rank $n$ is supersolvable if and only if
  it is $S_n$ \textrm{EL}-shellable.
\end{thm}

\subsection{Poset of forests}

\label{forests}

\begin{figure}
  \begin{center}
    \leavevmode 
    \epsfig{file=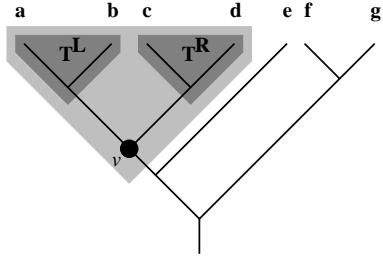,width=5cm} 
    \caption{The subtree $T_v$, and its parts $T^R$ and $T^L$.}
    \label{subtree}
  \end{center}
\end{figure}

A \textit{tree} is a leaf-labeled rooted binary tree and a
\textit{forest} is a set of such trees. Vertices are either inner
vertices (valence $3$) or leaves and roots (valence $1$). By
convention, edges are oriented towards the root. Leaves are
bijectively labeled by a finite set. Trees and forests are pictured
with their roots down and their leaves up, but are not to be
considered as planar. A leaf is an \textit{ancestor} of a vertex if
there is a path from the leaf to the root going through the vertex. If
$F_1,F_2,\dots,F_k$ are forests on $I_1,I_2,\dots,I_k$, let $F_1
\sqcup F_2 \sqcup \dots \sqcup F_k$ be their disjoint union. For a
forest $F$, we denote by $\V(F)$ the set of its {\em inner vertices}
and by $\LE(F)$ the set of {\em leaves}. The number of trees in a
forest $F$ on $I$ is the difference between the cardinal of $I$ and
the cardinal of $\V(F)$. By a {\em subtree} $T_v$ we mean the union of
all paths starting from any vertex $v$ and going up to the leaves.
Note that any subtree can be further divided in two parts denoted by
$T^L$ and $T^R$ as shown in Figure \ref{subtree}.

Following \cite{fred}, we introduce a partial order on the \textit{set
  of forests on $I$} denoted by $\For(I)$. 
\bigskip

\noindent{\bf Definition} 
Let $F$ and $G$ be forests
on the label set $I$. Then $F \leq G$ if there is a topological map
from $F$ to $G$ with the following properties:
\begin{enumerate}
\item[D1.] It is increasing with respect to orientation towards the root.
\item[D2.] It maps inner vertices to inner vertices injectively.
\item[D3.] It restricts to the identity of $I$ on leaves.
\item[D4.] Its restriction to each tree of $F$ is injective.
\end{enumerate}
In fact, such a topological map from $F$ to $G$ is determined up to
isotopy by the images of the inner vertices of $F$. One can recover
the map by joining the image of an inner vertex of $F$ in $G$ with
the leaves of $G$ which were its ancestor leaves in $F$.

\begin{figure}
  \begin{center}
    \leavevmode 
    \epsfig{file=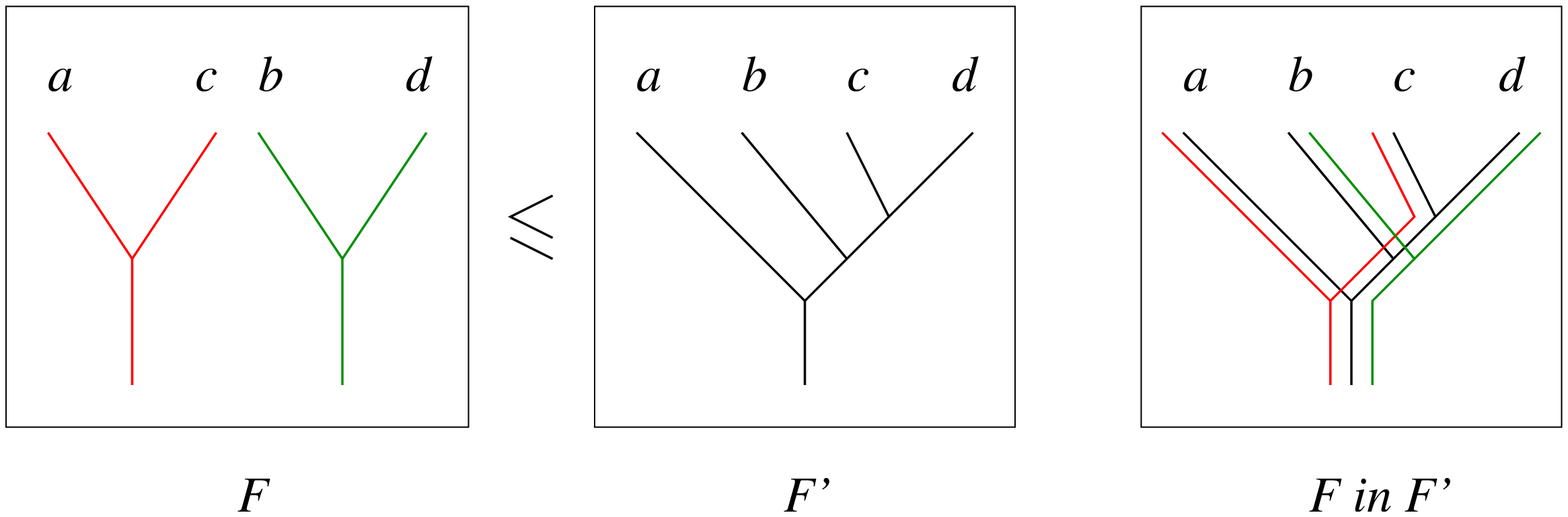,width=7cm} 
    \caption{$F \leq F^{\prime}$.}
    \label{exemple}
  \end{center}
\end{figure}

\begin{figure}
  \begin{center}
    \leavevmode 
    \epsfig{file=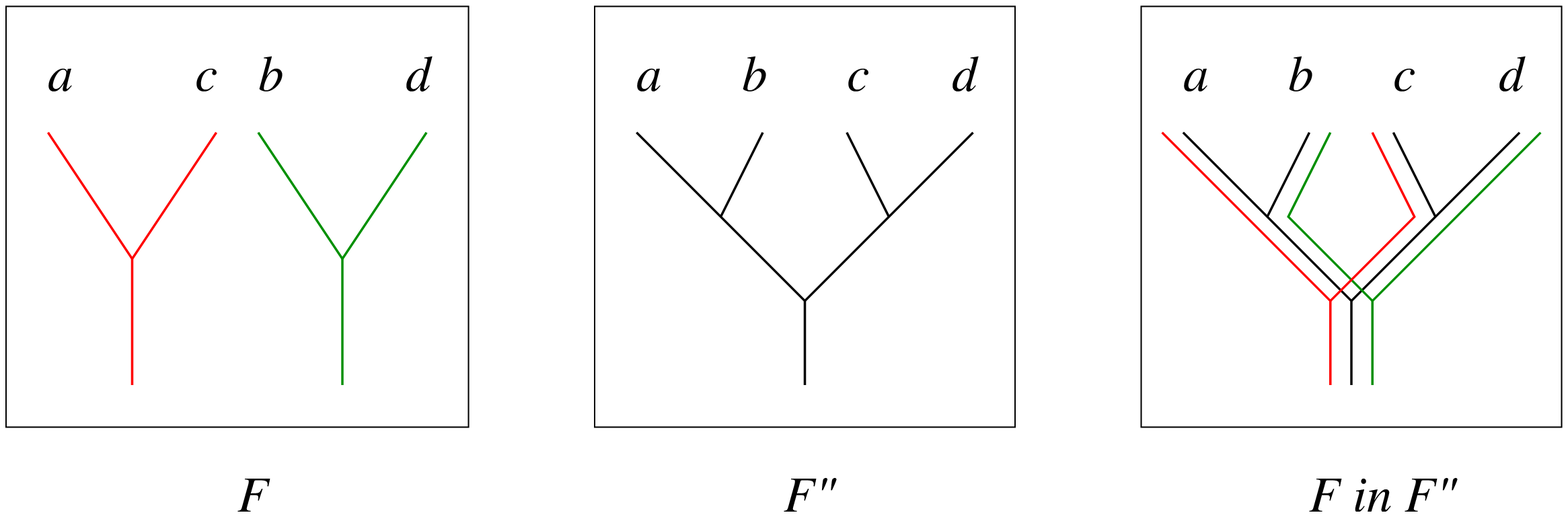,width=7cm} 
    \caption{$F$ and $F''$ are not comparable.}
    \label{contrex}
  \end{center}
\end{figure}

The following proposition can be found in \cite[Proposition 3.1]{fred}.
\begin{pro}
  \label{graded}
  The poset $\For(I)$ is graded by the number of inner vertices.
\end{pro}

It was proved in \cite{fred} that the maximal elements of the poset
$\For(I)$ are the trees. The forest without inner vertices is the
unique minimal element and is denoted by $\zero$. For any $J \subseteq I$, we denote by $|_J$ the tree such that $\V(|_J)=\emptyset$ and $\LE(|_J)=J$. Note that $\zero=|_I$.

\section{Intervals are lattices}

\label{latticity}

In this section we fix a finite set of leaves $I$ of cardinality $n+1$
and consider a tree $T$ on $I$. We study the interval $[\zero,T]$ that
is a graded bounded subposet of $\For(I)$. Our main goal is to show
that $[\zero,T]$ is a lattice.

\smallskip 

Any two distinct leaves $i,j \in I$ determine an inner vertex
$v_{(i,j)} \in \V(T)$, as the intersection of the two paths starting
from these leaves and going down to the root. Sometimes we will write
$i \stackrel{v}{\longleftrightarrow} j$ instead of $v=v_{(i,j)}$. For
any $J \subseteq I$, let
\begin{equation*}
  \SM(J):=
  \{v \in \V(T) \text{ : } v=v_{(i,j)} \text{ for some distinct } i,j \in J\}.
\end{equation*}

\begin{rem}\label{primo}
  For any subset $J \subseteq I$, it is easy to see that $|\SM(J)|=|J|-1$.
\end{rem}

\begin{lem}
  \label{sub}
  For any $J \subseteq I$, there exists a unique tree $T_J$ on $J$
  such that
  \[T_J\sqcup |_{I \setminus J} \leq T.\]
\end{lem}
\begin{proof}
  We define $T_J$ to be the union of all the paths starting from the
  leaves in $J$ and going down to the root. It is easy to check that
  all conditions in the definition of the partial order of forests are
  satisfied.
\end{proof}

\begin{rem}\label{arbre}
  Let $J_1 \subseteq J_2$ be two subsets of $I$. Then $T_{J_1} \sqcup |_{I \setminus J_1} \leq T_{J_2} \sqcup |_{I \setminus J_2}$
\end{rem}

The following definition is crucial in the rest of this paper.
\smallskip

Let $\pi=(\pi_1,\ldots,\pi_k)$ be a partition of $I$. We say that
$\pi$ is {\em $T$-admissible} if and only if
$\SM(\pi_i)\cap\SM(\pi_j)=\emptyset$ for all $i \neq j \in [k]$. We
denote the set of all $T$-admissible partitions of $I$ by
$\Ad(T)$.\\
For example, let $T=F''$ be the tree in Figure \ref{contrex} on
the set $I=\{a,b,c,d\}$. Then $\{\{a,b\},\{c,d\}\}\in \Ad(T)$, but
$\{\{a,c\},\{b,d\}\}$ is not a $T$-admissible partition of $I$, as in
fact $\SM(\{a,c\})=\SM(\{b,d\})=v_{(a,c)}$.

It is easy to see that $\Ad(T)$ is a poset by refinement order
$\leq_r$, \textit{i.e.} $(\pi_1, \ldots,\pi_n) \leq_r (\tau_1, \ldots,\tau_m)$
if and only if each block $\pi_i$ is contained in
some block $\tau_j$.\\
For example $\{\{a\},\{b,c\},\{d\}\} \leq_r \{\{a\},\{b,c,d\}\}$.

\smallskip

Let $F \in [\zero,T]$, $F=T_1\sqcup \ldots \sqcup T_k$, we define
\[\Pi(F):=(\pi_1,\ldots,\pi_k),\]
where $\pi_i:=\LE(T_i)$ for all $i \in [k]$.\\
Note that $\Pi(F)$ is a $T$-admissible partition by condition D2. 

\begin{pro}
  The map $\Pi:([\zero,T],\leq) \longrightarrow (\Ad(T),\leq_r)$ is an
  isomorphism of posets.
\end{pro}
\begin{proof}
  First we prove that $\Pi$ is a bijection. For every $\pi=(\pi_1,
  \ldots,\pi_k)\in \Ad(T)$, let
  \begin{equation}
    \label{gamma}
    \Gamma(\pi):=T_{\pi_1}\sqcup\ldots\sqcup T_{\pi_k},
  \end{equation}
  where
  each tree $T_{\pi_i}$ is defined by Lemma \ref{sub}.\\
  It is clear that $\Pi\circ \Gamma =\Id$. By the uniqueness in Lemma
  \ref{sub}, it follows that $\Gamma \circ \Pi = \Id$, and so $\Gamma$ is the inverse of $\Pi$.\\
  Now let $F,G \in [\zero,T]$ with $F \leq G$. Then, by condition
  D4, for all $T_F \in F$ there exists a $T_G
  \in G$ such that $\LE(T_F)\subseteq \LE(T_G)$. It follows that
  $\Pi(F)\leq_r \Pi(G)$. Conversely, if $\pi \leq_r \pi^{\prime}$,
  then, by Remark \ref{arbre}, we have $\Gamma(\pi) \leq
  \Gamma(\pi^{\prime})$. This concludes the proof.
\end{proof}
From now on, forests in $[\zero,T]$ and $T$-admissible partitions are
identified via the bijection $\Pi$.

\smallskip

We are ready to state and prove the main theorem of this section.

\begin{thm}
  For each tree $T$ on the set $I$, the interval $[\zero,T]$ is a
  lattice.
\end{thm}
\begin{proof}
  As the interval has a $\one$, by Proposition \ref{lattice} it
  suffices to prove that each $F,G \in [\zero,T]$ have a meet. Let
  $\Pi(F)=\pi=(\pi_1, \ldots,\pi_n)$ and $\Pi(G)=\tau=(\tau_1,
  \ldots,\tau_m)$. We show that the meet of $\pi$ and $\tau$ as
  partitions, defined by
  \[\pi\wedge\tau:=(\pi_1\cap \tau_1)\cup(\pi_1\cap
  \tau_2)\cup\ldots\cup (\pi_n\cap \tau_1)\cup\ldots \cup (\pi_n\cap
  \tau_m),\] is also in $\Ad(T)$. For every $(i,j) \neq
  (i^{\prime},j^{\prime})\in [n]\times [m]$ we have that \[\SM(\pi_i
  \cap \tau_j) \cap \SM(\pi_{i^{\prime}}\cap \tau_{j^{\prime}})
  \subseteq \SM(\pi_i)\cap \SM(\tau_j) \cap \SM(\pi_{i^{\prime}}) \cap
  \SM(\tau_{j^{\prime}})=\emptyset,\] because $\pi$ and $\tau$ are in
  $\Ad(T)$, hence either $\SM(\pi_i)\cap\SM(\pi_{i^{\prime}})$ or
  $\SM(\tau_j)\cap\SM(\tau_{j^{\prime}})$ is empty. It is immediate to
  see that $\pi \wedge \tau$ is the meet also in $\Ad(T)$, hence
  $\Ad(T)$ is a lattice and we are done.
\end{proof}

\section{$S_n$ EL-labelings on $[\zero,T]$}

\label{snelling}

In this section we introduce an edge-labeling on the poset $[\zero,T]$
and prove that it is an $S_n$ EL-labeling. By Theorem \ref{MN} it
follows that the lattice $[\zero,T]$ is supersolvable.

\smallskip
 
A partial order $\preceq$ is defined on the set $\V(T)$ in the
following way. 

\bigskip
\noindent{\bf Definition} A vertex $v$ is smaller than a vertex $v'$, denoted 
by $v \preceq v'$, if $v'$ is on the path between the root and $v$.
Any total order extending this partial order on $\V(T)$ is called a
\textit{nice} total order, still denoted by $\preceq$.

\begin{figure}
  \begin{center}
    \leavevmode 
    \epsfig{file=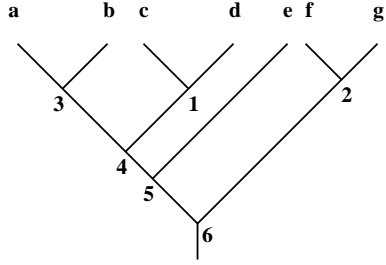,width=5cm} 
    \caption{Example of nice total order on $\V(T)$.}
    \label{totalorder}
  \end{center}
\end{figure}

\bigskip Using a nice total order, one can label the inner vertices by
integer numbers from $1$ to $n$. From now on, inner vertices and
labels are identified in this way using a fixed nice total order. Note
that the bottom vertex is the maximum element for the order $\preceq$.
An example is drawn in Figure \ref{totalorder}.

\bigskip

Now we introduce an edge-labeling as follows. First remark that for
all $F\leq G \in [\zero,T]$, one has $\V(F) \subseteq \V(G) \subseteq
\V(T)$. Moreover if $F \lhd G$, by Proposition \ref{graded}, there
exists a unique $v \in \V(G)$ such that $\V(G)=\V(F) \cup \{v\}$.
\bigskip

\noindent{\bf Definition}  Let $F \lhd G \in [\zero,T]$. Then we define $\la:\{(F,G)\; : \; F \lhd G\}
\rightarrow \NN$ by
\[\V(G)=\V(F) \cup \{\la(F,G)\},\]
where $\la(F,G)$ is the label of $v$.\\

\bigskip

\begin{figure}
  \begin{center}
    \leavevmode 
    \epsfig{file=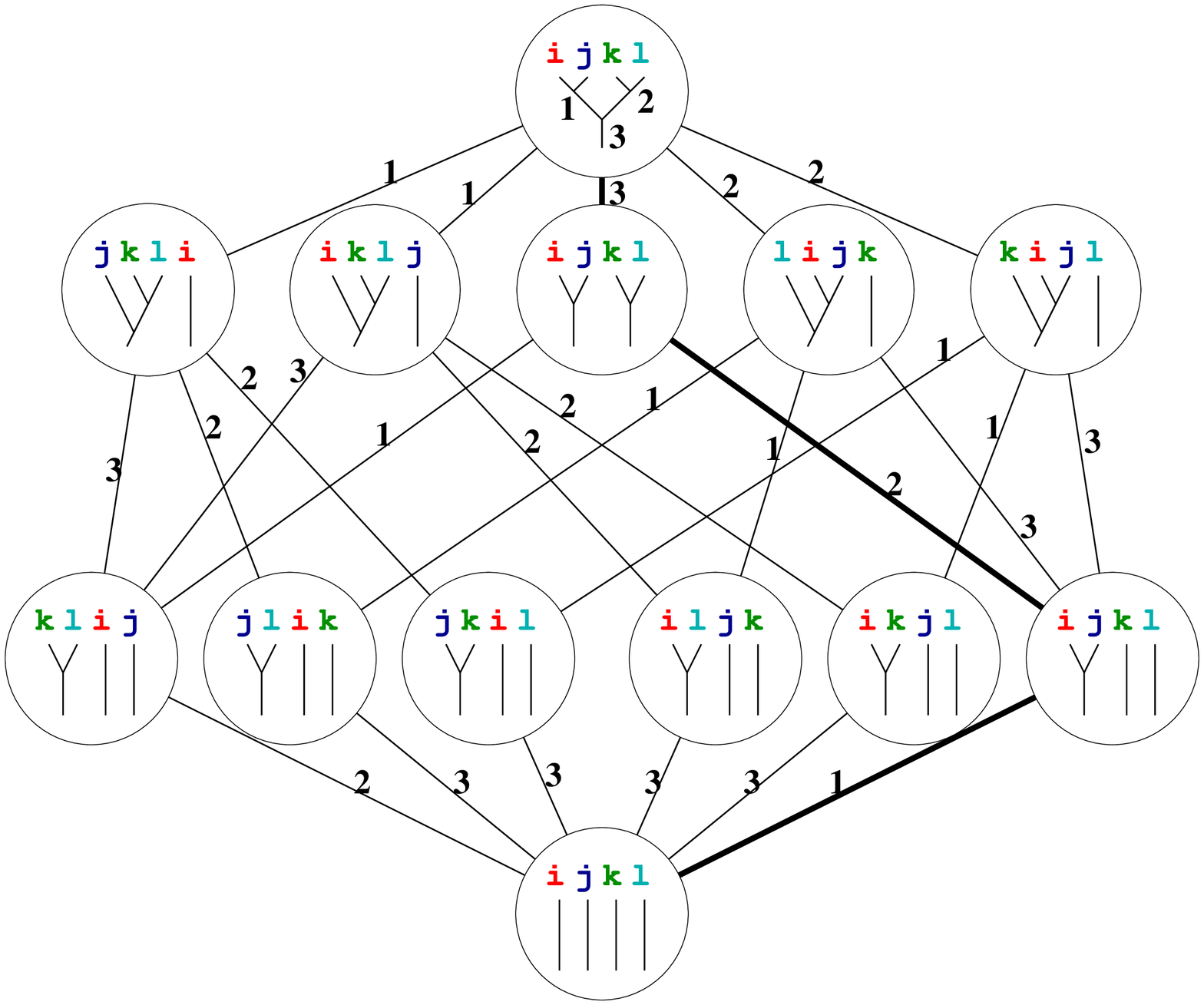,width=10.5cm} 
    \caption{$S_3$ EL-labeling of the interval $[\zero,T]$.}
    \label{shelling}
  \end{center}
\end{figure}

An example of this edge-labeling is shown in Figure \ref{shelling}.
The proof of the following lemma is immediate.
\begin{lem}
  \label{perm}
  The label of a maximal chain of $[F,G]$ is a permutation of the set
  $\V(G) \setminus \V(F)$.
\end{lem}

\begin{lem}
  \label{chain}
  For each $F \in [\zero,T]\setminus \{T\}$, there exists a unique $G \in
  [\zero,T]$ covering $F$ such that
  \[\la(F,G)={\rm min}(\V(T) \setminus \V(F)).\]
\end{lem}
\begin{proof}
  Let $\Pi(F)=\pi$ and let $v_0:={\rm min}(\V(T) \setminus \V(F))$.
  Consider the two subtrees starting from $v_0$, as explained in \S
  \ref{forests}, denoted $T_{v_0}^L$ and $T_{v_0}^R$. We show that
  $\LE(T_{v_0}^R)$ is contained
  in one part of $\pi$.\\
  Each $w \in \V(T_{v_0}^R)$ is such that $w \prec v_0$. It follows
  that $w \in \V(F)$ by minimality of $v_0$. Let $i \neq j \in
  \LE(T_{v_0}^R)$, then there is $v \in \V(T_{v_0}^R) \subseteq \V(F)$
  such that $i\stackrel{v}{\longleftrightarrow}j$. Hence $i,j$ are in
  the same part of $\pi$. Therefore $\LE(T_{v_0}^R)$ is contained in
  only one part of $\pi$ denoted by $\pi_R$. The same result is true
  for $T_{v_0}^L$, and we denote the corresponding part by $\pi_L$. As
  $v_0\not\in \V(F)$, the parts $\pi_L$ and $ \pi_R$ are distinct. We
  define a new partition
  \[\pi^{\prime}:=(\pi_L \sqcup \pi_R,\pi_1,\ldots,\pi_{k}),\]
  where $\pi_j$ are the remaining parts of $\pi$. From now on, we
  denote $\pi_L \sqcup \pi_R$ by $\pi_{LR}$.\\
  To show that
  $\pi^{\prime} \in \Ad(T)$, it suffices to prove that
  \begin{equation}
    \label{eq}
    \SM(\pi_{LR}) \cap \SM(\pi_j) = \emptyset, \;\;\; {\rm for \; all
      \;} j \in [k].
  \end{equation}
  We have that
  $\SM(\pi_{LR}) \supseteq \SM(\pi_L)\cup\SM(\pi_R)\cup\{v_0\}$. On the other hand, by Remark \ref{primo}, we have that $|\SM(\pi_L)|+|\SM(\pi_R)|+1=|\SM(\pi_{LR})|$, and so we have an equality. \\
  Now, for any $j \in [k]$, the vertex $v_0$ is not in $\SM(\pi_j)$,
  because all the ancestors of $v_0$ are in $\pi_L$ or in $\pi_R$.
  Hence condition $(\ref{eq})$ is verified and the proof of theorem
  follows by defining $G:=\Gamma(\pi_{LR},\pi_1,\ldots,\pi_{k})$,
  where $\Gamma$ is defined in (\ref{gamma}).
\end{proof}

The preceding lemma can be extended as follows.

\begin{pro}
  \label{chain2}
  For each $F, H \in [\zero,T]$, $F< H$ there exists a unique $G \in
  [\zero,T]$ covering $F$ such that
  \[\la(F,G)={\rm min}(\V(H) \setminus \V(F)).\]
\end{pro}
\begin{proof}
  If $H=T$ then the result is given by Lemma \ref{chain}. Otherwise
  let $H=H_1 \sqcup H_2\sqcup \ldots \sqcup H_k$, where $H_j$ is a
  tree for all $j\in [k]$. Since $F \leq H$, we have $F=F_1 \sqcup
  F_2\sqcup \ldots \sqcup F_k$ where $F_j$ is a forest, for all $j\in
  [k]$. It was observed in \cite[Proposition 2.1]{fred} that the
  interval $[F,H]$ is isomorphic to $\prod_{j=1}^k[F_j,H_j]$. Let
  $v_1:={\rm min}(\V(H) \setminus \V(F))$. We have $\V(H)=\V(H_1) \cup
  \V(H_2)\cup \ldots \cup \V(H_k)$ and, after re-ordering, we can
  assume that $v_1 \in \V(H_1)$. Then, by Lemma \ref{chain} applied to
  $[F_1,H_1]$, there exists a unique $G_1 \in [F_1,H_1]$ covering
  $F_1$ such that $\la(F_1,G_1)=v_1$. Define $G=G_1\sqcup F_2 \sqcup
  \ldots \sqcup F_k$ in $[F,H]$. Then $G$ is the unique forest of
  $[F,H]$, covering $F$, such that $\la(F,G)=v_1$. This concludes the
  proof.
\end{proof}

\begin{thm}
  \label{labeling}
  The lattice $[\zero,T]$ is \textrm{EL}-shellable. 
\end{thm}
\begin{proof}
  By Lemma \ref{perm}, for any interval $[F,G]$ of $[\zero,T]$, the
  unique possible increasing label for a saturated chain from $F$ to
  $G$ is given by the unique increasing permutation of the elements of
  $\V(G)\setminus \V(F)$.\\
  Then Proposition \ref{chain2} implies that there exists an unique
  chain $m$ from $F$ to $G$ with this label. The other maximal chains
  of $[F,G]$ are labeled by different permutations, which are
  lexicographically greater than the increasing one.\\
  Hence the edge-labeling $\lambda$ is an EL-labeling.
\end{proof}

\begin{cor}
  The lattice $[\zero,T]$ is supersolvable.
\end{cor}
\begin{proof}
  By Theorem \ref{labeling}, $\la$ is an EL-labeling and by Lemma
  \ref{perm}, $\la(m)$ is a permutation of $[n]$ for each maximal
  chain $m$. Hence $\lambda$ is an $S_n$ EL-labeling and the result
  follows from Theorem \ref{MN}.
\end{proof}
  
\begin{rem}
  \label{semi}
  Note that $[\zero,T]$ is not semimodular in general. For example,
  the atoms $\{\{j,k\},\{i\},\{l\}\}$ and $\{\{i,l\},\{j\},\{k\}\}$ in
  Figure \ref{shelling} do not satisfy the condition
  (\ref{semimodular}).
\end{rem}

\section{Characteristic polynomials}

In this section, we recover the results of \cite{fred} concerning the
characteristic polynomials of the intervals $[\zero,T]$. Note that, by
Remark \ref{semi}, the well-known theorem of Stanley \cite[Theorem
4.1]{staSS} (see also \cite[Theorem 6.2]{sagan99}) on the
factorization of the characteristic polynomials of semimodular
supersolvable lattices, does not apply. We use instead a stronger
theorem due to Blass and Sagan \cite{blassagan}.

\subsection{LL-lattices}

Recall that the characteristic polynomial of a graded finite lattice
$L$ of rank $n$ is
\begin{equation*}
  \chi_L(t)=\sum_{y \in L}\mu(\zero,y) t^{n-\rk(y)},
\end{equation*}
where $\mu$ is the Möbius function of $L$ and $\rk(y)$ is the rank of
$y$.

\smallskip

Following \cite{blassagan}, we define an element $x$ of a
lattice $L$ to be {\em left-modular} if, for all $y \leq z$,
\[ y \join (x \wedge z) = (y \join x) \wedge z. \]
A maximal chain $m \in {\mathcal M(L)}$ is said to be {\em
  left-modular} if all its elements are left-modular.
\begin{rem}
  \label{Mchain}
  From \cite[Proposition 2.2]{staSS}, it follows that if $L$ is a
  supersolvable lattice then its $M$-chain is left-modular.
\end{rem}

\smallskip

Any maximal chain $m: \zero =x_0 \lhd x_1 \lhd \cdots \lhd
x_n=\hat{1}$ defines a partition of the set of atoms $\AAA$ into
subsets called {\em levels} indexed by $i \in [n]$: 
\begin{equation*}
  \label{defia}
  \AAA_i=\{ a \in \AAA \text{ : } a \leq x_i \text{ and }a\not\leq x_{i-1}\}.
\end{equation*}

The partial order $\lhd_m$ on $\AAA$ {\em induced} by the maximal
chain $m$ is defined by
\[ a \lhd_m b \; {\rm if \; and \; only \; if } \;  a \in \AAA_i \; {\rm and}
\; b \in \AAA_j \; {\rm with} \; i <j.\] 
This partial order should not be confused with the covering
relation.

Then the following is called the {\em level condition} with respect to
$m$:
\[{\rm if } \; a_0\lhd_m a_1 \lhd_m a_2 \lhd_m \cdots \lhd_m a_k, \; {\rm then} \; a_0 \not \leq \bigvee_{i=1}^{k} a_i.\] 

A lattice $L$ having a maximal chain $m$ that is left-modular and
satisfies the level condition is called an LL-\textit{lattice}.

The following theorem is due to Blass and Sagan \cite[Theorem 6.5]{blassagan}.
\begin{thm}
  \label{blass}
  Let $P$ be an {LL}-lattice of rank $n$. Let $\AAA_i$ be the
  levels with respect to the left-modular chain of $P$. Then
  \begin{equation*}
    \chi_P(t)=\prod_{i=1}^{n} (t-| \AAA_i |).
  \end{equation*}
\end{thm}

\subsection{Factorization of characteristic polynomials}

A tree $T$ with $n$ inner vertices and leaf set $I$ is fixed. A nice
total order on $\V(T)$ is chosen, defining an edge-labeling as in
\S\ref{snelling}.

\smallskip

The set $\AAA$ of atoms of $[\zero,T]$ is the set of pairs $(i,j)$ of
distinct elements of $I$. To each atom $(i,j)$ is associated an inner
vertex $v_{(i,j)}$ of $T$ as defined in \S \ref{latticity}. The covering
edge $\zero \lhd (i,j)$ is labeled by the integer in $[n]$
corresponding to $v_{(i,j)}$ in the chosen total order on $\V(T)$.

\begin{pro}
  \label{join}
  Let $a_1,a_2,\dots,a_k \in \AAA$ with pairwise distinct
  vertices $v_1,v_2,\dots,v_k$ in $\V(T)$. Then $\V(a_1 \join a_2 \join \dots \join a_k)=\{v_1,v_2,\dots,v_k\}$.
\end{pro}
\begin{proof}
  Let $V=\{v_1,v_2,\dots,v_k\}$. Let
  $\pi^{(1)},\pi^{(2)},\dots,\pi ^{(k)}$ be the partitions of $I$
  associated to $a_1,a_2,\dots,a_k$. Let $\pi$ be the join $\pi ^{(1)}
  \join \pi ^{(2)} \join \dots \join \pi ^{(k)}$ in the lattice of
  partitions. We want to show
  that $\pi \in \Ad(T)$ and that $\V(\pi)=V$.\\
  Let $p$ be a part of $\pi$. Let $V_p$ be the set of vertices in $V$
  whose corresponding atoms in $\{a_1,\dots,a_k\}$ have their leaves in $p$. Observe that the sets $V_p$
  form a partition of $V$ because atoms in $\{a_1,\dots,a_k\}$ have
  pairwise distinct vertices. Let $v$ be a vertex in $\SM(p)$. This
  means that there exists $i,j$ in $p$ such that
  $i\stackrel{v}{\longleftrightarrow}j$. As $p$ is a part of a join,
  there exists a chain
  \begin{equation*}
    i=i_0\stackrel{t_0}{\longleftrightarrow}i_1
    \stackrel{t_1}{\longleftrightarrow}i_2\dots
    i_{\ell-1}\stackrel{t_{\ell-1}}{\longleftrightarrow}
    i_\ell \stackrel{t_\ell}{\longleftrightarrow} i_{\ell+1}=j,
  \end{equation*}
  where each $i_r \stackrel{t_r}{\longleftrightarrow} i_{r+1}$ is an
  atom in $\{a_1,\dots,a_k\}$ with vertex in $V_p$.\\
  In the rest of the proof, the symbol $\preceq$ stands for the
  partial order introduced in \S \ref{snelling}.\\
  Let us prove by induction on the length $\ell$ of the chain that
  there exists $\theta_\ell$ in $V_p$ such that $\theta_\ell \succeq
  t_0$ and
  $\theta_\ell \succeq t_\ell$.\\
  If $\ell=0$, then one can take $\theta_0=t_0$. Assume that there
  exists $\theta_{\ell-1}$ in $V_p$ such that $\theta_{\ell-1} \succeq
  t_0$ and $\theta_{\ell-1} \succeq t_{\ell-1}$. The path joining the
  leaf $i_\ell$ to the root contains the vertices
  $t_{\ell-1}$,$t_\ell$ and hence also by induction hypothesis the
  vertex $\theta_{\ell-1}$. Either $t_{\ell} \preceq \theta_{\ell-1}$,
  and one can take $\theta_{\ell}=\theta_{\ell-1}$ or $t_{\ell}
  \succeq \theta_{\ell-1}$ and one can take $\theta_{\ell}=t_{\ell}$.
  This concludes the induction.\\
  Therefore $\theta_{\ell} \in V_p$ is such that
  $i\stackrel{\theta_\ell}{\longleftrightarrow}j$. Hence
  $\theta_\ell=v \in V_p$ and so $\SM(p) \subseteq V_p$. The
  converse inclusion is clear.\\
  Now let $p$ and $p'$ be two different parts of $\pi$. Then
  $\SM(p)\cap \SM(p')=V_p \cap V_{p'}$ is empty. Hence $\pi$ is
  $T$-admissible.\\ We have proved that $\pi$ is $T$-admissible and
  that the vertices of $\pi$ are exactly $V$. It follows that $\pi$
  defines the join $a_1 \join \dots \join a_k$ in $[\zero,T]$ and the
  proposition is proved.
\end{proof}

Define another partition of $\AAA$ indexed by $i \in [n]$:
\begin{equation*}
  \BBB_i=\{ a \in \AAA \text{ : } \lambda(\zero,a)=i \}.
\end{equation*}

Let $m:\zero=x_0\lhd x_1 \lhd \dots \lhd x_n=T$ be the fixed modular chain of
$[\zero,T]$, \textit{i.e.} the unique increasing maximal chain for the fixed
labeling.

\begin{lem}
  \label{atomic}
  Let $i \in [n]$. For each $j \in [i]$, let $a_j$ be an atom in
  $\BBB_j$. Then
  \begin{equation*}
    x_i=a_1 \join a_2 \join \dots \join a_i.
  \end{equation*}
\end{lem}
\begin{proof}
  The proof is by induction on $i$. By Proposition \ref{chain2},
  $x_1=a_1$ is the unique atom in $\BBB_1$. Assume that $x_{i-1}=a_1
  \join \dots \join a_{i-1}$. Then $a_1 \join \dots \join a_{i-1}
  \join a_{i}$ is $x_{i-1} \join a_{i}$ and has rank $i$ by
  Proposition \ref{join}. Moreover we have that
  $\lambda(x_{i-1},x_{i-1} \join a_{i})=i$. By uniqueness in
  Proposition \ref{chain2}, it follows that $x_{i}=x_{i-1} \join
  a_{i}$.
\end{proof}

\begin{lem}
  \label{aabb}
  Let $\AAA_i$ be the levels with respect to $m$. Then for each $i\in
  [n]$,
  \begin{equation*}
    {\AAA}_i={\BBB}_i.
  \end{equation*}
\end{lem}
\begin{proof}
  It suffices to prove that
  \begin{equation*}
     \{ a \in \AAA \text{ : } a \leq x_i \}=\{ a \in \AAA \text{ : }
    \lambda(\zero,a)\in [i] \}.  
  \end{equation*}
  If $a \leq x_i$, then $\lambda(\zero,a)$ is one of the vertices of
  $x_i$, \textit{i.e.} belongs to $[i]$. Conversely, take any atom $a$ with
  $\lambda(\zero,a)$ in $[i]$. Choose other atoms to have one atom in
  each $\BBB_j$ for $j \in [i]$. Then, by Lemma \ref{atomic}, $x_i$ is
  the join of $a$ and the other chosen atoms, so $a \leq x_i$.
\end{proof}

\begin{pro}
  \label{lllattice}
  The lattice $[\zero,T]$ is an \textrm{LL}-lattice.
\end{pro}
\begin{proof}
  This lattice is supersolvable, so by Remark \ref{Mchain} the
  $M$-chain is a left-modular chain. It remains to check the level
  condition. Take atoms $a_0,a_1,\dots, a_k$ which belongs to pairwise
  different $\AAA_i$. By Lemma \ref{aabb}, these atoms belong to
  pairwise different $\BBB_i$. Then by Proposition \ref{join} the set
  of vertices of the join $a_1 \join \dots \join a_k$ does not contain
  the vertex of the atom $a_0$. This ensures the level condition.
\end{proof}

Now we are ready to state and prove the main result of this section,
which was already proved in \cite[Theorem 4.6]{fred}.

\begin{figure}
  \begin{center}
    \leavevmode 
    \epsfig{file=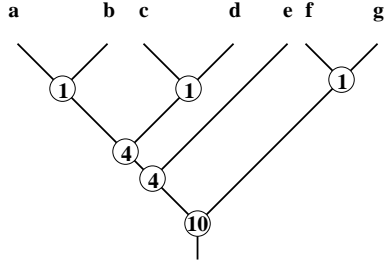,width=5cm} 
    \caption{Example of roots of the characteristic polynomial.}
    \label{arbrexpo}
  \end{center}
\end{figure}

\begin{thm}
  The characteristic polynomial of $[\zero,T]$ is
  \begin{equation*}
    \chi_{[\zero,T]}(t)=\prod_{v \in \V(T)} (t-e(v)),
  \end{equation*}
  where $e(v)$ is the product of the number of left ancestor leaves of
  $v$ by the number of right ancestor leaves of $v$.
\end{thm}
\begin{proof}
  By Proposition \ref{lllattice}, one can apply Theorem \ref{blass} to
  $[\zero,T]$. Let us count the number of elements of $\AAA_i$ for
  each $i$. By Lemma \ref{aabb}, this is equal to the cardinality of
  $\BBB_i$. Let $v$ be the vertex of $T$ with index $i$. It is easy to
  see that the number of atoms in $\BBB_i$ is the number of left
  ancestor leaves of $v$ times the number of right ancestor leaves of
  $v$.
\end{proof}
For example, the characteristic polynomial of the interval $[\zero,
T]$ where $T$ is the tree in Figure \ref{arbrexpo} is
$\chi_{[\zero,T]}(t)=(t-1)^3(t-4)^2(t-10)$.
\bibliographystyle{plain}
\bibliography{tree.bib}

\end{document}